\documentclass[12pt,twoside,reqno,a4paper]{amsart}
\usepackage[cp1251]{inputenc}
\usepackage[russian,english]{babel}
\usepackage{amsmath,amssymb,graphicx,color}

\usepackage[colorlinks,pdfauthor={Yu. Brezhnev},
            pdfwindowui=false,
            bookmarksopen=false,bookmarks=false]{hyperref}

\def~{\unskip\nobreak\hspace{0.27778em}\ignorespaces}
\def\,{\ifmmode\mskip+1.5mu\else\kern+.08333em\fi\relax}
\def\!{\ifmmode\mskip-1.5mu\else\kern-.08333em\fi\relax}
\newdimen{\FontSize} \delimiterfactor=990
\def\cdot{\,{\mathchar"2201}\,}

\makeatletter

\def\dddot#1{{\mathop{{}#1}\limits^{\vbox to-1.5\ex@{\kern-\tw@\ex@
\hbox{\larger[2]\rm.\kern-0.12em.\kern-0.12em.}\vss}}}} \makeatother

\def\ds{\displaystyle}

\def\s{\scriptscriptstyle}
\def\ie{i.\,\,e.}

\def\re{{\,\mathrm{e}}}
\def\ri{{\mathrm{i}}}

\def\wpp {\wp\Smaller{'}}

\def\ep {e\Smaller{'}}
\def\epp{e\Smaller{'\mskip-1.2mu'}}

\def\II{{\textbf{\textsc{I\kern -0.15em I}}}}
\def\III{{\textbf{\textsc{I\kern -0.17em I\kern -0.17em I}}}}
\def\IIIab{\III\!\raise0.2ex\hbox{$\begin{smallmatrix}
\alpha\\[-0.1ex]\varrho\end{smallmatrix}$}\kern-0.1em}

\def\DEF{\mathrel{\vcenter{\hbox{$:$}}{=}}}

\def\=={\mathrel{\phantom{=}}}

\def\sm{\mathrel{\vcenter{\hbox{\scalebox{0.8}[1]{$\scriptstyle-$}}}}}

\def\?{\textrm{\protect\footnotesize$\red\mathchar"446$}}

\def\vphi{\smash[b]{\raise-0.22ex\hbox{\Smaller[2]{\s\boldsymbol-}}
  \mkern-8.00mu\raise0.435ex\hbox{\scalebox{1}[0.9]{$\varphi$}}}\relax}
\def\vpsi{\smash[b]{\raise-0.22ex\hbox{\Smaller[2]{\s\boldsymbol-}}
  \mkern-8.95mu\raise0.408ex\hbox{\scalebox{1}[0.9]{$\psi   $}}}\relax}

\def\russian{\selectlanguage{russian}}
\def\english{\selectlanguage{english}}
\def\END{\end{document}}

\def\Mfrac#1#2{\def\FRAC{\hbox{\smaller[1]$\ds\frac{#1}{#2}$}}
 \ifdim\FontSize=10pt\raise0.057ex\FRAC\else
{\ifdim\FontSize=11pt\raise0.050ex\FRAC\else\raise0.051ex\FRAC\fi}\fi}

\def\mfrac#1#2{\def\FRAC{\hbox{\smaller[2]$\ds\frac{#1}{#2}$}}
 \ifdim\FontSize=10pt\raise0.1166ex\FRAC\else
{\ifdim\FontSize=11pt\raise0.104ex\FRAC\else\raise0.0968ex\FRAC\fi}\fi}

\newcommand{\AB}[2]{\nobreak\mskip-1mu
\raise0.12ex\hbox{\smaller[2]\renewcommand{\arraystretch}{0.42}%
\setlength{\arraycolsep}{0pt}\raise0.16ex\hbox{$\big[$}%
$\begin{array}{c}\scriptstyle#1\\ \scriptstyle#2\end{array}$%
\raise0.16ex\hbox{$\big]$}}\mskip-1mu\relax}

\newcommand{\mbig}[2][4]{\vcenter{\hbox{%
\ifcase#1$#2$\or\small$\big#2$\or$\big#2$\or\large$\big#2$%
\or\footnotesize$\Big#2$\or\small$\Big#2$\or$\Big#2$\or\large$\Big#2$%
\or$\bigg#2$\or\large$\bigg#2$\or$\Bigg#2$\or\large$\Bigg#2$%
\or\Large$\Bigg#2$\else\huge$\bigg#2$\fi}}\relax}

\newdimen{\upp} \newdimen{\Dim} \newdimen{\DDim}
\newcommand{\pow}[4][\displaystyle]{
\settowidth{\upp}{\hbox{$#1 {}^{#4}$}}%
\settowidth{\Dim}{\hbox{$#1{#2}$}}%
\settowidth{\DDim}{\hbox{$#1{#2}_{#3}$}}%
\addtolength{\DDim}{-\Dim} {#1{#2}_{#3}\kern-\DDim {}^{#4}}
\addtolength{\DDim}{-\upp}\ifdim\DDim>0pt \kern+\DDim \fi\relax}

\makeatletter 
\providecommand{\larger}[1][1]{
\ifmmode\relax\else\@tempcnta\ifx\@currsize\normalsize 4\else
 \ifx\@currsize\small 3\else\ifx\@currsize\footnotesize 2\else
 \ifx\@currsize\large 5\else\ifx\@currsize\Large 6\else
 \ifx\@currsize\LARGE 7\else\ifx\@currsize\scriptsize 1\else
 \ifx\@currsize\tiny  0\else\ifx\@currsize\huge 8\else
 \ifx\@currsize\Huge  9\else 4 \fi\fi\fi\fi\fi\fi\fi\fi\fi\fi
\advance\@tempcnta#1\relax\ifnum\@tempcnta<\z@ \@tempcnta\z@ \fi
\ifcase\@tempcnta\tiny\or\scriptsize\or\footnotesize\or\small
 \or\normalsize\or\large\or\Large\or\LARGE\or\huge\else\Huge\fi\fi}
\makeatother

\providecommand{\smaller}[1][1]{\larger[-#1]}

\newcommand{\Smaller}[2][1]{\hbox{\smaller[#1]$#2$}}

\makeatletter \def\@dotsep{.8} \makeatother 

\def\red{\color{red}}
\def\blue{\color{blue}}

\newcommand{\Width}[2][\displaystyle]
{\settowidth\Dim{\hbox{$#1#2$}}\texttt{\tiny\blue$($\the\Dim$)$}}
\newcommand{\Height}[2][\displaystyle]
{\settoheight\Dim{\hbox{$#1#2$}}\texttt{\tiny\blue$($\the\Dim$)$}}
\newcommand{\Depth}[2][\displaystyle]
{\settodepth\Dim{\hbox{$#1#2$}}\texttt{\tiny\blue$($\the\Dim$)$}}
\usepackage{hyperref}

\newcommand{\textcenter}[3][Yu.~Brezhnev]{
\pagestyle{myheadings}
\markboth{\hfill\textsc{#1}\hfill}{\hfill\textsc{#1}\hfill}
\ifdim#2=0mm\textwidth=\textwidth\else\textwidth=#2\fi
\oddsidemargin=\paperwidth%
\addtolength{\oddsidemargin}{-\textwidth}%
\addtolength{\oddsidemargin}{-2in} \oddsidemargin=0.5\oddsidemargin
\evensidemargin=\oddsidemargin
\ifdim#3=0mm\textheight=\textheight\else\textheight=#3\fi
\topmargin=\paperheight%
\addtolength{\topmargin}{-\textheight}%
\addtolength{\topmargin}{-2in}%
\addtolength{\topmargin}{-\headsep}%
\addtolength{\topmargin}{-\headheight}%
\addtolength{\topmargin}{-\footskip}%
\topmargin=0.5\topmargin}

 \textcenter{150mm}{237mm} 

\theoremstyle{plain}

\theoremstyle{remark}

\allowdisplaybreaks[4] 

\def\uu{{\mathfrak{u}}}

\def\bz{{\boldsymbol{z}}}

\def\II{{\boldsymbol{\mathrm{I\kern -0.15em I}}}}
\def\III{{\boldsymbol{\mathrm{I\kern -0.18em I\kern -0.18em I}}}}
\def\A{\mbox{\small$A$}}
\def\B{\mbox{\small$B$}}
\def\sA{\sqrt{\!\!\A\,}}
\def\sB{\sqrt{\!\B\,}}
\def\kpm{k^{\s\pm}}

\begin{document}

\author{\sc Yu.~V.~Brezhnev}

\title[Uniformization of Abelian
integrals]{On uniformizable representation\\for Abelian
integrals}

\date{\today}

\thanks{Research supported by the Federal Targeted Program
under contract 02.740.11.0238}


\begin{abstract}
We show how the unformizable representation for Abelian integrals of
nontrivial genera arise. The technique  makes use of the famous
Chudnovsky's Fuchsian linear differential equations and their relation
to the sixth Painlev\'e transcendent.
\end{abstract}

\hfill {\smaller[2]{\bf Painlev\'e equations and Related Topics,}
199--208. {\smaller \copyright} De~Gruyter 2012}\bigskip

\maketitle\thispagestyle{empty}


\section{Introduction}
In this note we exhibit first examples of uniformizable
$\tau$-representation for functions having an additive automorphic
property, \ie, Abelian integrals on 1-dimensional orbifolds of a
negative constant curvature. Orbifolds are generalizations of a notion
Riemann surface for the case when fundamental group of this `manifold'
has elements of finite or formally infinite order. These are called
usually conic singularities or punctures. Matrix representations for
fundamental groups of 1-dimensional orbifolds are described by Fuchsian
linear ordinary  differential equations (ODEs) of 2nd order. Their
singularities are precisely the conic points on these orbifolds.

In work \cite{br1} we showed how Fuchsian  equations
\begin{equation}\label{P}
p\,\Psi''+p'\,\Psi'+(x+A)\,\Psi=0\,,\qquad p\DEF
x\,(x-\alpha)(x-\beta)
\end{equation}
arise in the theory of Painlev\'e equations and in note \cite{br2} we
announced a first explicit example of a solvable allied Schwarz ODE
$$
[\uu,\tau]=-2\,\wp(2\,\uu)\,,\qquad \uu(\tau)=
\frac{\vartheta_3(\tau)}{\vartheta_2(\tau)}\!\cdot\!
{}_2F_1\!\!\!\left(\Mfrac12,\Mfrac14;\Mfrac54\Big|
\Mfrac{\vartheta_3^4(\tau)}{\vartheta_2^4(\tau)}\right),
$$
where $\uu$ is an elliptic holomorphic integral. Below are some details
on these  examples and methods of getting explicit formulae.

\section{Schwarz equation and equations on tori}
As is well known Schwarz's equations are the 3rd order ODEs coming from
linear ODEs of the form $\Psi_{\mathit{xx}}=\frac{1}{2}\,Q(x)\,\Psi$.
The ratio of its two linearly independent solutions
$$
\tau=\frac{\Psi_1(x)}{\Psi_2(x)}
$$
defines $\tau$ as function of $x$ and conversely. We may write down an
autonomic ODE defining $x$ as function of $\tau$. If, as usual,
$$
\{f,z\}\DEF
\frac{f_{\!\mathit{zzz}}}{f_z}-
\frac32\frac{\pow{f}{\mathit{zz}}{2}}{\pow{f}{z}{2}}
$$
defines  the standard Schwarz derivative then the inverted function
$x=\chi(\tau)$ satisfies the ODE
\begin{equation}\label{S}
[x,\tau]= Q(x)\,,
\end{equation}
wherein the following notation has been adopted
$[x,\tau]\DEF-\{\tau,x\}$, that is
$$
[x,\tau]\DEF \frac{x_{\tau\tau\tau}}{\pow{x}{\tau}{3}}-\frac{3}{2}\,
\frac{\pow{x}{\tau\tau}{2}}{\pow{x}{\tau}{4}}\,.
$$

We shall deal with the two examples of Fuchsian equations of the class
\eqref{P}:
\begin{alignat*}{3}
x\,(x-1)(x+1)\,&\Psi''+{}&(3\,x^2-1)\,&\Psi'+{}&(x+0)\,\Psi=0\,,\\
x\,(x^2+3\,x+3)\,&\Psi''+{}&(3\,x^2+6\,x+3)\,&\Psi'+{}&(x+1)\,\Psi=0\,,
\end{alignat*}
They belong to the set of four equations known as Chudnovsky's ones
\cite{chud}. With the help of well-known linear transformation
$\Psi\mapsto \psi$ we can transform these equations into the normal
form
\begin{alignat}{2}
\psi''&=-\frac14\,\frac{(x^2+1)^2}{x^2(x-1)^2(x+1)^2}\,
\psi\,,\label{1'}\\[0.3em]
\psi''&=-\frac14\,\frac{(x+1)(x+3)(x^2+3)}{x^2(x^2+3\,x+3)^2}\,
\psi\,.\label{2'}
\end{alignat}

It is known that each of these equations  has the four parabolic
singularities on Riemann sphere and, thereby, define two punctured
spheres; the simplest orbifolds of genus zero. Halphen used original
trick to relate these equations  with elliptic functions corresponding
to algebraic curves of the form $y^2=x\,(x-1)(x+1)$ and
$y^2=x\,(x^2+3\,x+3)$. They have the standard Weierstrassian models
$$
y^2=4\,x^3-4\,x\quad\text{and}\quad
y^2=4\,x^3-4
$$
respectively, \ie, Gauss' lemniscate and the equi-anharmonic elliptic
curve.

Let us make the substitutions $x=\wp(\uu;g_2^{},g_3^{})$ in equations
\eqref{1'} and \eqref{2'}, where invariants $g_2^{}$, $g_3^{}$ are
chosen according to the Weierstrassian models above. We shall obtain
Fuchsian equations in variable $\uu$ and, then, Schwarz's equation of
the form \eqref{S}. An easy computation gives
\begin{equation}\label{tori}
[\uu,\tau]=-2\,\wp(2\,\uu;g_2^{},g_3^{})\,.
\end{equation}
and this equation is equivalent to equations on tori considered for the
first time in classical work \cite{keen}.  Thus, if we find solutions
$\uu=\uu(\tau)$ to this equation we obtain nontrivial examples of
uniformizable representation for the object $\uu$ which is an Abelian
integral, because
$$
\uu=\wp^{\sm1}(x;g_2^{},g_3^{})\,.
$$
On the other hand, above mentioned Fuchsian and Schwarz's equations
correspond to orbifolds with punctures and therefore their curvature is
not equal to zero but is a negative constant. The theory of such
orbifolds is very nontrivial. Automorphic functions on them (Klein's
Hauptmoduln) are known, a few as they are, but examples of additively
automorphic objects are absent hitherto.

\section{Holomorphic elliptic integrals and hypergeometric functions}
\subsection{Lemniscate}
We have in this case
\begin{equation*}\label{ppp}
\pm\uu=\int\limits_\infty^{\;x}\!\!\!
\frac{du}{\sqrt{4\,u^3-4\,u\,}}=\cdots
\end{equation*}
It follows that
\begin{equation*}\label{ind}
\cdots =\frac12\!\int\limits_\infty^{\;x}\!\!
u^{\sm\frac12}_{\mathstrut}(u^2-1)^{\sm\frac12}_{\mathstrut}\,du
\end{equation*}
which is  a particular case of a consequence of the integral definition
to the hyper\-geo\-metric ${}_2F_1$-function. More precisely, the
standard integral definition of the ${}_2F_1$-functions is through the
definite integral \cite[Sect.~2.1.3]{bateman2}
$$
{}_2F_1(\alpha,\beta;\gamma|z)\DEF\frac{\Gamma(\gamma)}
{\Gamma(\beta)\,\Gamma(\gamma-\beta)}\,\int\limits_0^{\;1}\!\!
u^{\beta\sm1}(1-u)^
{\gamma\sm\beta\sm1}(1-z\,u)^{\sm\alpha}\, du
$$
but in some cases this definition may be rewritten in terms of
non-definite integrals. Changing the integration variable $u\mapsto
z\,u$, the upper limit $u=1$ gets mapped to $u=z$. Putting further
$\gamma-\beta-1=0$, we arrive at a definition of this particular case
for the ${}_2F_1$-function through  the \textit{indefinite} integral.
Clearly, all the other cases are in fact variations of this
scheme\footnote{To all appearances, the variable upper limit integral
formulae for definition to the ${}_2F_1$-functions was observed for the
first time in the 1876 dissertation by Tikhomandritski\u\i\
\cite[p.~78]{tih} before the known 1881 Goursat dissertation.}. We
obtain (see also tables in \cite{tables})
\begin{equation}\label{Fab}
\int\limits_0^{\;z}\!\!u^{\alpha\sm1}(u-1)^{\sm\beta}\,du=
\frac{\re^{\pi\ri\beta}}{\alpha}\,z^\alpha\cdot
{}_2F_1 (\beta,\alpha;\alpha+1|z)\,,\qquad \mathrm{Re}\,(\alpha)>0\,,
\end{equation}
and, therefore,
\begin{equation}\label{table}
\int\limits_\infty^{\;z}\!\!
u^{\alpha\sm1}(u-1)^{\sm\beta}\,du=
\frac{z^{\alpha-\beta}}{\alpha-\beta}\cdot
{}_2F_1 (\beta,\beta-\alpha;\beta-\alpha+1|z^{\sm1})\,,\qquad
\mathrm{Re}\,(\beta-\alpha)>0\,.
\end{equation}
The lemniscate integral under question is thus transformed into
\begin{equation}\label{pp}
\uu=\frac{1}{\sqrt{x\,}}
\cdot{}_2F_1\!
\mbig[6](\Mfrac12,\Mfrac14;\Mfrac54\Big|\Mfrac{1}{x^2}\mbig[6]),
\end{equation}
We also know that $\tau$-representation for all the Chudnovsky
equations are expressed through the standard elliptic modular
functions, namely, Jacobi's $\vartheta$-constants. By this means,
correlating these two points, we obtain an explicit solution to the
Schwarz equation \eqref{tori}.

Let us use the Hauptmodul $x=\chi(\tau)$ for the 1st Chudnovsky
equation \eqref{1'}:
$$
\chi(\tau)=\frac{\vartheta_2^2(\tau)}{\vartheta_3^2(\tau)}\,.
$$
where
$$
\vartheta_2(\tau)\DEF
\re^{\frac14\pi\ri\tau}_{\mathstrut}
\sideset{}{_k}\sum\limits_{-\infty}^{\infty}\!\!
\re^{(k^2+k)\pi\ri\,\tau}\,,\quad
\vartheta_3(\tau)\DEF \sideset{}{_k}\sum\limits_{-\infty}^{\infty}\!\!
\re^{k^2\pi\ri\,\tau}\,,\quad
\vartheta_4(\tau)\DEF
\sideset{}{_k}\sum\limits_{-\infty}^{\infty}{}\!\! (-1)^k\,
\re^{k^2\pi\ri\,\tau}\,.
$$
To put it differently  this function solves equation
$$
[x,\tau]=-\frac12\,\frac{(x^2+1)^2}{(x^3-x)^2}\,.
$$
Substituting this $\chi(\tau)$ into \eqref{pp}, we get that function
$$
\uu(\tau)=\frac{\vartheta_3(\tau)}{\vartheta_2(\tau)}\!\cdot\!
{}_2F_1\!\!\!\left(\Mfrac12,\Mfrac14;\Mfrac54\mbig[7]|
\Mfrac{\vartheta_3^4(\tau)}{\vartheta_2^4(\tau)}\right)
$$
solves Eq.~\eqref{tori} under $(g_2^{},g_3^{})=(4,0)$.

\subsection{Equi-anharmonic curve}
First we shift $x$-variable $x=z-1$ in Eq.~\eqref{2'} to obtain the
canonical form $y^2=4\,z^3-4$ with $(g_2^{},g_3^{})=(0,4)$. Halphen's
transformation and Hauptmodul $x=\chi(\tau)$ \cite{maier} in this case
have the form
$$
z=\wp(\uu;0,4)\,,\qquad z=9\,\frac{\eta^3(9\,\tau)}{\eta^3(\tau)}+1\,,
$$
where
$\eta(\tau)\DEF\re^{\frac{\pi\ri}{12}\tau}_{\mathstrut}\prod_k(1-\re^{2\pi\ri
k\tau})$ is the Dedekind eta-function. This  Hauptmodul $z(\tau)$
satisfies the equation
$$
[z,\tau]=-\frac12\,\frac{z\,(z^3+8)}{(z^3-1)^2}
$$
and one can show that the change $z\mapsto \uu$ above transforms this
Schwarz equation into Eq.~\eqref{tori}. We have, according to
\eqref{table},
$$
\pm\uu=\int\limits_\infty^{\;z}\!\!\!
\frac{du}{\sqrt{4\,u^3-4\,}}=\frac{1}{\sqrt{z}}\cdot
{}_2F_1\!
\mbig[6](\Mfrac12,\Mfrac16;\Mfrac76\Big|\Mfrac{1}{z^3}\mbig[6])
$$
and, hence,
$$
\uu=\bigg(9\,\frac{\eta^3(9\,\tau)}{\eta^3(\tau)}+1\bigg)
^{\!\!\!\!\sm\frac12}\cdot
{}_2F_1\!
\mbig[8](\Mfrac12,\Mfrac16;\Mfrac76\mbig[7]|
\mbig[7]\{9\,\Mfrac{\eta^3(9\,\tau)}{\eta^3(\tau)}+1\mbig[7]\}^
{\!\!\!\sm3}\mbig[8])\,.
$$
This is a very nontrivial exercise to check directly that this function
solves Eq.~\eqref{tori}.

\textit{Remark}. All the Chudnovsky Hauptmoduln are single-valued
functions and, hence, in spite of a square root in the last formula,
this expression provides a single-valued object in the neighborhood of
point $z=\infty$. Indeed, making the change $z^2=\bz^{\sm1}$ in
Schwarz's equation for $z$ above, we get
$$
[\bz,\tau]=-\frac{1}{2}\frac{1}{\bz^2}+\cdots
$$
and, hence, $\bz=\bz(\tau)$ has an exponentially  single-valued
behavior in $\tau$ about point $z=\infty$ and ${}_2F_1$-function is, by
definition, a single valued Taylor series at the origin $\bz=0$.

We can, however, get an explicit root-free solution to this problem
with the help of the  following trick. Use  formula \eqref{table} and
the fact that holomorphic integral is defined up to an additive
constant:
$$
\uu=\int\limits_\infty^{\;x}\!\!\!
\frac{du}{\sqrt{4\,u^3-4\,}}=\int\limits_\infty^{\;0}\!\!\!
\frac{du}{\sqrt{4\,u^3-4\,}}+\int\limits_0^{\;x}\!\!\!
\frac{du}{\sqrt{4\,u^3-4\,}}\,.
$$
The first of these integrals is a transcendental constant $\uu_{\s0}$
such that $\wp(\uu_{\s0})=0$, \ie,  zero of the Weierstrass
$\wp(z;0,4)$-function. It is computed into elliptic integrals
$$
\begin{aligned}
\uu_\mathrm{o}&=\frac{\ri}{2\sqrt[4]{3}}\,F\Big(\sqrt[4]{3}\,
\big(\sqrt{3}-1\big);
\mfrac14\,\big(\sqrt{6}+\sqrt{2}\big)\Big)-\frac{1}{\sqrt[4]{3}}\,
K\Big(\mfrac{\sqrt{3}-1}{2\,\sqrt{2}}\Big)=\\
&=\frac{\ri}{6}\,\mathrm{B}\Big(\Mfrac16,\Mfrac13\Big)=
\ri\cdot1.402\,182\,105\,325\ldots,
\end{aligned}
$$
where $F$ and $K$ are the standard non-complete and complete elliptic
integrals \cite{bateman2} and B is the Euler beta-function. The second
integral, upon application of \eqref{Fab}, becomes
$$
\int\limits_0^{\;x}\!\!\! \frac{du}{\sqrt{4\,u^3-4\,}}=
\frac{\ri}{2}\,x\cdot
{}_2F_1\!
\mbig[6](\Mfrac12,\Mfrac13;\Mfrac43\Big|x^3\mbig[6])\,.
$$
Turning to the Hauptmodul (variable) $z(\tau)$, we get finally
$$
\pm\uu(\tau)=\uu_\mathrm{o}+\frac{\ri}{2}\,
\bigg(9\,\frac{\eta^3(9\,\tau)}{\eta^3(\tau)}+1\bigg)
\cdot
{}_2F_1\!
\mbig[8](\Mfrac12,\Mfrac13;\Mfrac43\mbig[7]|
\mbig[7]\{9\,\Mfrac{\eta^3(9\,\tau)}{\eta^3(\tau)}+1\mbig[7]\}^
{\!\!3}\mbig[8])\,.
$$
Single-valuedness of this expression is now obvious; it solves
Eq.~\eqref{tori}.

\section{Abelian integrals for genus $g>1$}

If $\alpha$, $\beta$ in formula \eqref{Fab} are the rational numbers
then one has an Abelian integral belonging to rationality that may have
genus greater than unity. The integral may be holomorphic, meromorphic,
or logarithmic one. If, on the other hand, we have a single-valued
$\tau$-representation for $z^\alpha$, then we get automatically an
explicit  $\tau$-representation for this integral. This is a nontrivial
consequence because, in general, Abelian integrals are not expressible
in terms of any known functions and, thereby, we get analogs of
Weierstrass' uniformization theory for elliptic curves.

Recall that Weierstrass' theory, in its complete description, includes
closed and self-contained collection of the base functions and
differentials
\begin{equation}\label{tor1}
w^2=4\,z^3-a\,z-b\,,\qquad
\big\{z=\wp(\tau)\,,\;w=\wpp(\tau)\big\}\,,\qquad
dz=\wpp(\tau)\,d\tau\,,
\end{equation}
and the two principal integrals
\begin{equation}\label{tor2}
\II\DEF\!\!\int\limits^{\,\,\,z}\!z\,\frac{dz}{w}=-\zeta(\uu)\,,\qquad
\III\DEF\frac12\!\int\limits^{\,\,\,z}\!\frac{w+\wpp(\alpha)}
{z-\wp(\alpha)}\,\frac{dz}{w}
=\ln\frac{\sigma(\uu-\alpha)}{\sigma(\uu)}+
\zeta(\alpha)\,\uu\,,
\end{equation}
as functions of the fundamental holomorphic object
\begin{equation}\label{tor3}
\uu\DEF\int\limits^{\,\,\,z}\!\frac{dz}{w}\,,\qquad\uu(\tau)=\tau\,.
\end{equation}
By this means, having the base holomorphic and meromorphic integrals we
can manipulate with this integrals and functions in exactly the same
manner as we do with Weierstrassian objects $\sigma$, $\zeta$, $\wp$,
and $\wpp$. Insomuch as no one explicit formula for such a theory was
known hitherto, it is, perhaps, not without interest to exhibit
examples.

\subsection{Higher genera. Examples}
Consider the set of integrals
$$
\int\limits_0^{\;z}\!\!
\frac{u^m\,du}{\sqrt[\uproot{2}\leftroot{-1}k]{u\,(u^4-1)}}\,.
$$
Looking them at the integrals for the algebraic irrationality
$w^k=z\,(z^4-1)$, we find that they can be either holomorphic or
meromorphic ones. For example, in case $k=2$ (hyperelliptic curves) we
have two base holomorphic integrals when $m=0$, $1$ and two base
meromorphic ones when $m=2$ or $3$. The only thing we need now is the
formula for Hauptmodul $z=z(\tau)$. It is known explicitly in cases
when $z$, as a solution to corresponding Fuchs--Schwarz's equation, has
only parabolic singularities; that is Fuchsian equation has punctures
at all the points $z=\{0,\pm1,\pm\ri,\infty\}$:
$$
z=\frac{\vartheta_2(\tau)}{\vartheta_3(\tau)}
$$
We derive from \eqref{Fab}
\begin{equation}\label{aa}
\begin{aligned}
\int\limits_0^{\;z}\!\!u^{\alpha\sm1}(u^n-1)^{\sm\beta}\,du
&=
\frac{\re^{\pi\ri\beta}}{\alpha}\,z^\alpha\cdot
{}_2F_1\mbig[6](\beta,\Mfrac{\alpha}{n};\Mfrac{\alpha}{n}+1
\mbig[5]|z^n\mbig[6])\,,\\
\int\limits_\infty^{\;z}\!\!u^{\alpha\sm1}(u^n-1)^{\sm\beta}\,du
&=
\frac{z^{\alpha-n\beta}}{\alpha-n\,\beta}\cdot
{}_2F_1\mbig[6](\beta,\beta-\Mfrac{\alpha}{n};\beta-\Mfrac{\alpha}{n}+1
\mbig[5]|\Mfrac{1}{z^n}\mbig[6])
\end{aligned}
\end{equation}
and $\alpha=m+1-\frac1k$. The value $n$ does not affect on
single-valuedness of the $\tau$-representation and under $k=2$ we have
to do the square root of $z(\tau)$. One can show that ratio of any two
Jacobi's $\vartheta$-constant is a complete square. In particular
$$
\sqrt{\frac{\vartheta_2(\tau)}{\vartheta_3(\tau)}}=\sqrt{2}\,
\frac{\vartheta_2(\tau)}{\vartheta_2\big(\frac{\tau}{2} \big)}\,.
$$
It follows that for the famous algebraic curve $w^2=z^5-z$ we can
obtain a complete set of base Abelian integrals.

\textbf{Theorem.} \textit{Every Abelian $($homolomorphic, meromorphic,
or logarithmic\/$)$ integral belonging to the algebraic irrationality
$w^2=z^5-z$ has a uniformizing $\tau$-representation through the
Jacobi's $\theta$-functions of the two base holomorphic integrals}.

\medskip
\noindent \textit{Proof}. This curve admits a representation in
form of two isomorphic  elliptic curves. Indeed, it is suffice to use
the well-known formulae by Jacobi--Legendre for covering of the curve
$$
y^2=x\,(x-1)(x-\A)(x-\B)(x-\A\B)\,.
$$
We  have in this case the following cover of torus $(\uu)$:
$$
\left\{
\begin{aligned}
\wp(\uu)&=-\frac{(\sA\pm\sB)^2}{(x-\A)(x-\B)}\,x-
\frac13\,(\kpm+1)\\[0.8em]
\wpp(\uu)&=2\frac{(\sA\pm\sB)^2}{\sqrt{(1-\A)(1-\B)}}\,
\frac{x\pm\sqrt{\!\!\A\B}}{(x-\A)^2(x-\B)^2}\,y
\end{aligned}
\right.,
$$
where
$$
\kpm=-\frac{(\sA\pm\sB)^2}{(1-\A)(1-\B)}\,.
$$
Then expression
$\wpp(\uu)^2=4\,\wp^3(\uu)-g_2^{}\,\wp(\uu)-g_3^{}$ is equivalent to
the hyperelliptic curve above.  Reduction of holomorphic differentials
has the form
$$
d\uu=\frac12\,\sqrt{(1-\A)(1-\B)}\!\cdot\!
(x\mp\sA\sB)\,
\frac{dx}{y}\,.
$$
and the curve $w^2=z^5-z$ corresponds to the following parameters
$$
\A=-1\,,\qquad\B=\ri\,,\qquad\kpm=\frac12\,(1\pm\sqrt{2})\,,
$$
$$
\begin{aligned}
\wpp{}^2&=4\,\wp^3-\frac53\,\wp\pm\frac{7}{27}\,\sqrt{2}=\\[0.5em]
&=4\Big\{\wp+\underbrace{\tfrac16(3\pm\sqrt{2})}_{\ds -e}\!\Big\}
\Big\{\wp-\underbrace{\big(\!{\pm}\tfrac13\sqrt{2}\big)}_
{\ds \ep}\!\Big\}
\Big\{\wp-\underbrace{\tfrac16(3\mp \sqrt{2})}_{\ds \epp}\!\Big\}\,.
\end{aligned}
$$
With use of \eqref{aa} we derive
$$
\mathfrak{U}=\int\limits_0^{\;z}\!\!
\frac{u^m\,du}{\sqrt{u\,(u^4-1)}}=\frac{2\,\sqrt{2}\,\ri}{2\,m+1}\,
\frac{\vartheta_2^{m+1}(\tau)}{\vartheta_3^m(\tau)\,
\vartheta_2\big(\frac{\tau}{2} \big)}\cdot
{}_2F_1\!
\mbig[8](\Mfrac12,\Mfrac14m+\Mfrac18;\Mfrac14m+\Mfrac98\mbig[7]|
\Mfrac{\vartheta_2^4(\tau)}{\vartheta_3^4(\tau)}\mbig[8])
$$
and these expressions provide not only two holomorphic integrals
$\mathfrak{U}_1$, $\mathfrak{U}_2$ ($m=0,1$) but the meromorphic ones
as well. Insomuch as we have holomorphic integrals as functions of
$\tau$ we can compute any other Abelian integral in terms of
2-dimensional $\Theta$-functions. On the other hand, by virtue of the
cover above, all the $\Theta$-functions are reducible to Jacobi's
$\theta$. It follows that all the Abelian integrals are expressed
through such $\theta$-functions having their arguments the two
holomorphic integrals above. These $\theta$-formulae can be transformed
into the Weierstrass $\wp$, $\wpp$, $\zeta$, $\sigma$ because all the
integrals are reducible to the objects  \eqref{tor1}--\eqref{tor3}.
\hfill $\blacksquare$

\end{document}